\documentclass[12pt]{article}
\usepackage[a4paper,margin=2.5cm]{geometry}
\usepackage[active]{srcltx}
\usepackage{geometry}

\usepackage{amsfonts}
\usepackage{amssymb}
\usepackage{amsfonts}
\usepackage{amsthm}
\usepackage{bbm}
\usepackage{bm}
\usepackage{amsmath}
\usepackage[normalem]{ulem}

\usepackage{longtable}

\usepackage{times}
\usepackage{fullpage}
\usepackage{lineno}

\usepackage{color,cite,graphicx}
\fussy
\usepackage{chapterbib}

\def\[{\begin{equation}}
\def\]{\end{equation}}
\newcounter{numer1}

\usepackage{multirow}
\usepackage{setspace}

\usepackage{subcaption}
\usepackage{graphicx}
\usepackage{amssymb}
\usepackage[all]{xy}

\usepackage{epsfig}

\usepackage{xcolor}

\usepackage{color}
\usepackage{hyperref}
\hypersetup{
    colorlinks=true, 
    linkcolor=blue, 
    urlcolor=red, 
    linktoc=all 
}

\usepackage{todonotes}

\newtheorem{theorem}{Theorem}
\newtheorem{lemma}{Lemma}
\newtheorem{proposition}{Proposition}
\newtheorem{coro}{Corollary}
\newtheorem{definition}{Definition}

\newtheorem{OBS}{Remark}

\newcommand{\im}{\mbox{\rm Im}}

\usepackage{bm}

\newcommand{\be}{\begin{equation}}
\newcommand{\ee}{\end{equation}}

\newcommand{\incorrect}[1]{{\color{black}#1}}
\newcommand{\suggestion}[1]{{\color{black}#1}}

\usepackage{color}

\title{Recovering sparse networks: Basis adaptation and stability under extensions}
\author{Marcel Novaes$^{1,2}$, Edmilson Roque dos Santos$^{1}$ and Tiago Pereira$^{1, 3}$ \\ {\small $^1$ Instituto de Ci\^encias Matem\'aticas e Computa\c{c}\~ao, Universidade de S\~ao Paulo, Brasil }\\
\small{$^2$ Instituto de F\'{i}sica, Universidade Federal de Uberl\^andia, Uberl\^andia, MG} \\
\small{$^3$ Department of Mathematics, Imperial College London, London SW7 2AZ, UK}} 

\begin{document}
\maketitle

\begin{abstract}
We consider the problem of recovering equations of motion from multivariate time series of oscillators interacting on sparse networks. We reconstruct the network from an initial guess which can include expert knowledge about the system such as main motifs and hubs. When sparsity is taken into account the number of data points needed is drastically reduced when compared to the least-squares recovery. We show that the sparse solution is stable under basis extensions, that is, once the correct network topology is obtained,  the result does not change if further motifs are considered. \end{abstract}


\section{Introduction}

Networks of interacting self-sustained oscillators have become a rich interdisciplinary topic, with applications ranging from neuroscience to physics and sociology \cite{RMP}. Across diverse applications the properties of the network may vary significantly. For example, the number of participants  ranges from a few to hundreds of thousands, 
and the interaction structure can consist of everyone interacting with everyone, or exhibit small-world properties, or be based on hierarchical structures among the participants \cite{newman2010}.

Once the mathematical description of the system is given, recent work has combined the theory of dynamical systems with graph theory to understand the impact of the network structure in the overall behavior.  This approach has been able to successfully demonstrate that the network structure can have systematic influences on properties such as synchronization \cite{eroglu2017,JEMS}.

In experiments, it is often impossible to directly determine the network structure, though. In fact, typically one has access to certain states of individual elements of the network, thus obtaining multivariate time series. A fundamental challenge is to recover the network interaction structure from data. This question has attracted much attention~\cite{timme2007, Domenico_2008,rosenblum_2017,kralemann2014reconstructing,pikovsky1,pikovsky2,pik2018,Eroglu_PRX_2020}.

Usually, the recovery uses prior expert knowledge \suggestion{of} possible network structures. From these guesses, one may extend the recovery reconstructing further interactions. This set of examples contains many important applications such that in neuroscience and engineering. 

In this work, we study the reconstruction of sparse networks. We start from a network seed that gives an approximation of the network to be recovered and extend the search for further connections. \suggestion{We show that $(i)$ by adapting the recovery to the dynamics, $(ii)$, the basis extension does not lead to prediction instability. We discuss the least square techniques are unstable under basis extension. A heuristic upshot of our study is that if the network is sparse and has $k\ll N$ links, where $N$ being the number of nodes in the network, then using sparse recovery we need only $O(N k)$ data points as opposed to least-square where we need $O(N^2)$.}

We will focus on the case when isolated dynamics of the nodes have a stable periodic motion and the interaction is weak. This is an interesting case, as the phase itself is not observed and thus we need to preprocess the data.

\section{Dynamics near a Hopf Bifurcation}

We consider the isolated dynamics of each node in the network to be near a Hopf-Andronov bifurcation, modelled by the Stuart-Landau equation
\begin{linenomath*}
\be \label{maineq}\dot{z_i}=F(z_i)=(1+{\rm j} \omega_i)z_i-|z_i|^2z_i,
\ee
\end{linenomath*}
where $z_i$ is a complex number.
Each isolated oscillator has an exponentially attractive periodic orbit with amplitude $1$ and frequency $\omega_i$ for $i = 1,\dots, N$.  The effect of a linear pairwise interactions is modelled as
\begin{linenomath*}
\be 
\dot{z_i}=F(z_i)+\alpha \sum_{k=1}^N C_{ik}(z_k-z_i)
\ee
\end{linenomath*}
for $i=1,\dots,N.$ Here, $\alpha$ denotes the coupling strength, assumed small. The connectivity matrix $C$ describes the interaction structure: $C_{ik}$ is $1$ if node $i$ is influenced by node $k$ and is $0$ otherwise. Notice that in the absence of linear terms, if nonlinear terms are included in the coupling this could lead to higher order resonances. However, we will consider only linear coupling which is enough to show how the recovery method works.

\subsection{Phase Dynamics} 

By introducing polar coordinates $z_i = r_i e^{{  \rm j} \theta_i}$ we can obtain the dynamics of amplitudes $r_i$ and phases $\theta_i$. As $\alpha$ is small, the network effect on the amplitudes is small, in fact, $r_i(t) = 1+ O(\alpha)$. The relevant dynamics generated by the network is encoded in the phases. The coupled phase equations \suggestion{to} leading order in $\alpha$ read as
\begin{linenomath*}
\be\label{phase1}
\dot{\theta_i}=\omega_i+\alpha \sum_{k=1}^n C_{ik} \sin(\theta_k-\theta_i).
\ee
\end{linenomath*}

\noindent {\bf Extracting phase from data.}  
In applications we do not have direct access to  $\theta_i(t)$, and may need to \suggestion{infer} another phase variable from a time series.  Let $x_i$ and $y_i$ denote, respectively, the real and imaginary parts of $z_i$ and that we assume that we only measure $x_i(t)$ for each oscillator. Thus, we have a multivariate time series for the network. To extract the phase from each time series we use the standard Hilbert transform
\begin{linenomath*}
\be
H(x_i(t)) =\frac{1}{\pi}  \mbox{p.v.} \int_{-\infty}^{+\infty}\frac{x_i(\tau)}{t-\tau}d\tau.
\ee
\end{linenomath*}
Thus using the analytic signal 
\begin{linenomath*}
\be s_i(t) = x_i(t) + {  \rm j}H(x_i(t)) = R_i(t) e^{{  \rm j} \vartheta_i(t)}
\ee
\end{linenomath*}
we can extract a phase $\vartheta_i(t)$ corresponding to the signal $x_i(t)$. \suggestion{Although} this phase is a surrogate and not necessarily equal to $\theta_i(t)$, meaningful dynamical information can be obtained from it. Once we have the phases $\vartheta_i$, their time derivatives are obtained numerically and a smoothing filter is applied to remove noise introduced in this process. 

\section{The recovery method}

\subsection{The basis functions}
The idea is to express the time derivatives of the phases, obtained from data, as linear combinations of certain functions. Here as we deal \suggestion{with} phases we use Fourier modes depending on all variables and on the differences of all variables,
\begin{linenomath*}
\begin{equation}\label{eq:phase_model}
    \dot{\vartheta}_i=\omega_i+\sum_\ell g^{(i)}(\vartheta_\ell)+\sum_{k,m} h^{(i)}(\vartheta_k,\vartheta_m),
\end{equation}
\end{linenomath*}
where 
\begin{linenomath*}
\begin{equation}\label{eq:isolated_function}
    g^{(i)}(\vartheta_\ell)=a_\ell^{(i)}\cos(\vartheta_\ell)+b_\ell^{(i)}\sin(\vartheta_\ell)
\end{equation}
\end{linenomath*}
is the isolated component and the coupling function is 
\begin{eqnarray} 
h^{(i)}(\vartheta_k,\vartheta_m)=c_{k,m}^{(i)}\cos(\vartheta_k-\vartheta_m)+d_{k,m}^{(i)}\sin(\vartheta_k-\vartheta_m),\quad k<m.
\end{eqnarray} 
The choice of coupling function $h$, depending only on phase differences, is motivated by the theory of phase reduction.  

The aim is to find the coefficients $\{a,b,c,d\}$ that provide a good approximation to the data $\dot{\vartheta}_i$.  We have $N$ time series for our $\vartheta_i$ variables, with $n$ points each, obtained with a fixed known sampling rate. With this data we form time-series for the $m = 1 + 2N + N(N-1)/2$ Fourier modes and arrange them as columns of a $n\times m$ matrix, we denote it as ${\Theta}$, so
\begin{linenomath*}
\be 
{\Theta} =
{\small \frac{1}{\sqrt{n}}
\left(
\begin{array}{ccccccc}
1 	   & \sin(\vartheta_{1}(t_1)) & \cdots & \sin(\vartheta_{N}(t_1))& \cos(\vartheta_{1}(t_1))&\cdots &\cos(\vartheta_{N-1}(t_1)-\vartheta_{N}({t_{1}})) \\
1      & \sin(\vartheta_{1}(t_2)) & \cdots & \sin(\vartheta_{N}(t_2))& \cos(\vartheta_{1}(t_2))&\cdots &\cos(\vartheta_{N-1}(t_2)-\vartheta_{N}({t_{2}})) \\
1      & \sin(\vartheta_{1}(t_3)) & \cdots & \sin(\vartheta_{N}(t_3))& \cos(\vartheta_{1}(t_3))&\cdots &\cos(\vartheta_{N-1}(t_3)-\vartheta_{N}({t_{3}})) \\
\vdots & \vdots                 & \vdots & \vdots                & \vdots & \ddots &\vdots         \\
1      & \sin(\vartheta_{1}(t_n)))& \cdots & \sin(\vartheta_{N}(t_n))& \cos(\vartheta_{1}(t_n))&\cdots &\cos(\vartheta_{N-1}(t_n)-\vartheta_{N}({t_{n}})) \\
\end{array}
\right)}. \nonumber
\ee
\end{linenomath*}
The problem of recovering the equations of motion can be formulated as the search for a $m\times n$ matrix of coefficients ${W}$ such that the equation 
\begin{linenomath*}
\be\label{sparse1}
{\Theta W}={V}
\ee
\end{linenomath*}
is satisfied, where  
\begin{linenomath*}
\be 
{V}=\left(
\begin{array}{ccc}
\dot \vartheta_1(t_1) & \cdots &  \dot\vartheta_N(t_1)	    \\
\vdots    &\ddots & \vdots\\
\dot \vartheta_1(t_n)  &   \cdots &  \dot\vartheta_N(t_n)	
\end{array}
\right)
\ee
\end{linenomath*}
is a $n\times N$ matrix of time series of derivatives, we apply smoothening to the derivatives.

The matrix $\Theta$ has all possible connections and because the network is sparse only a subset will contribute. We will denote by 
\begin{itemize}
\item[] -- $A$ a subset of columns of $\Theta$ that contain the expert guess. 
\item[] -- $B$  further columns we wish to probe. 
\end{itemize}
Without loss of generality (up to relabelling nodes) we assume that $A$ 
correspond to the first $p$ columns of $\Theta$. Next we consider the \suggestion{concatenation} of $[A,B]$ of the matrices $A$ and $B$ and 
consider the problem 
\begin{linenomath*}
\begin{equation}
[A,B] w = v \nonumber
\end{equation}
\end{linenomath*}
\noindent
where $v$ is one of the columns of the matrix $V$.
The vector of coefficients $w$ can be decomposed 
in terms of the action of $A$ and $B$
\begin{linenomath*}
\begin{equation}
w= \left(
\begin{array}{c}
x \\
y
\end{array}
\right).\nonumber
\end{equation}
\end{linenomath*}
The remaining exposition will address two problems: How to find the vector coefficients $x$, and the effect of the basis 
extension $B$ on the solution $x$. 

\subsection{The minimization}

Consider the problem of finding the 
vector of coefficients starting from the expert guess
\begin{linenomath*}
\begin{equation}
A x = v \nonumber
\end{equation}
\end{linenomath*}
The least squares approximation  provides the vector ${x}$ that minimizes the $L_2$ error
\begin{linenomath*}
\begin{equation}
\min_{x\in \mathbb{R}^{p}} \| A {x}-{v}\|_2. \nonumber
\end{equation}
\end{linenomath*}
A major advantage of this $L_2$ minimization is that the unique solution has a closed form,
 \begin{linenomath*}
\begin{equation}
\label{eq:L2solution}
{x}_0=A^+ {v} 
\end{equation}
\end{linenomath*}
where $A^+=(A^{\dagger} A)^{-1}A^{\dagger}$ is the pseudoinverse of $A$ and $\dagger$ denotes the transpose.  

Kraleman et al.\cite{pikovsky1,pikovsky2} have  used $L_2$ minimization to recover the topology of networks with up to nine oscillators. For a brief review, see \cite{pik2018}. Notice that, although this approach minimizes the euclidean error, it may not be an optimal solution with respect to other criteria, specially when the smallest singular value of $A$ becomes small. To obtain a well conditioned matrix the size of the time series needs to be significantly large. 

Denote $\im~A$ the image of the matrix $A$. Let us consider the case $n > p$, if $v \in \im ~A$ the system of equations has a unique solution and it is independent of the minimization. As the data is subjected to fluctuations\suggestion{,} in general 
\begin{linenomath*}
\begin{equation}
v = b + z \nonumber
\end{equation}
\end{linenomath*}
where $b \in \im ~A$ and $z \in (\im ~A)^{\perp}$, the orthogonal complement, with $\| z\|_2 \le \varepsilon$, for some small $\varepsilon > 0$ capturing the fact that fluctuations are small. 

\subsection{Finding Sparse Solutions}
For example, we may want a sparse solution, i.e. a vector $x$ with a few non-zero elements. This will indeed be the case when the network has sparse connectivity (such as the star network we shall consider, which has only $N$ connections out of a total of $N(N-1)$ possibilities). \par 
Sparsity can be measured in terms of the condition where
\begin{linenomath*} 
\be
\| x \|_0 = \mbox{number of nonzero elements of } x
\ee
\end{linenomath*}
should be as small as possible. Finding a sparse solution is a combinatorial NP-hard problem and not tractable. \suggestion{When the matrix $\Theta$ has some additional structure, namely it satisfies the restricted isometry property (RIP) \cite{candes1}}, it is well known that a valid heuristics to obtain sparse solutions is to include in the minimization process a penalization on the $L_1$ 
norm,
\begin{linenomath*}
\be
\| x \|_1=\sum_{i=1}^m |x_i|
\ee
\end{linenomath*}
\suggestion{and} consider 
\begin{linenomath*}
\begin{equation}\label{P1}
\underset{\tilde{w} \in \mathbb{R}^m}{\mbox{min}} \| \tilde{w} \|_1 \mbox{~ subject to ~}  \| \Theta \tilde{w}- v \|_2<\varepsilon,
\end{equation}
\end{linenomath*}
for some small $\varepsilon$, where we are still considering $\Theta = [A,B]$. This is known as basis pursuit denoising \cite{donoho}. The solution to this problem can be obtained by quadratic programming. This is the idea behind the Matlab package ``{\rm l1magic}''\footnote{https://statweb.stanford.edu/~candes/l1magic/}. However, there is a small technical drawback here, which is that to start the search for a minimal solution one needs a seed, and this is usually the $L_2$ solution similar to Equation~(\ref{eq:L2solution}). In situations when this $L_2$ solution is a poor choice (see at Section~\ref{sec:L2_unstable}), the algorithm may not be successful (and finding other clever seeds is a challenging problem). 

Another approach is the LASSO algorithm (least absolute shrinkage and selection operator), which we shall adopt. It works by computing solutions to
\begin{linenomath*}
\be\label{lasso}
\min_{w \in \mathbb{R}^m} \| \Theta w- v \|_2^2+\lambda \|w\|_1 
\ee
\end{linenomath*}
for a series of values of $\lambda$. When $\lambda$ is large, the solution approaches the null vector. When $\lambda$ is gradually decreased, each previous solution is a good seed for a new minimization process that finds sparse solutions. If $\lambda$ becomes too small, sparsity is no longer promoted. 

Intermediate values of $\lambda$ therefore lead to solutions that come close to minimizing 
$\| \Theta x- v\|_2$, while at the same time being significantly sparse. The actual value of $\lambda$ is selected by a process of $k$-fold cross validation, in which: the data is split into $k$ equal-sized parts; a solution is found using all but the $l$th part; a prediction error is computed when predicting the behavior on the $l$th part; the errors are added for $1\le l\le k$ to form the total prediction error; the value of $\lambda$ is chosen to minimize the total prediction error.  
Later we will see that by our Theorem \ref{thm:coherence} once we establish an adapted basis,  LASSO is not affected by the poor conditioning of $\Theta$ and performs significantly better when data acquisition time is short.

\section{Numerical experiments}
\subsection{Results for a directed star}

We consider a directed star motif for a paradigm. It consists of a central node driving to $N-1$ peripheral nodes, as shown in Figure~\ref{fig:CompLASSOL2}. Since every node's dynamics is only influenced by node 1, the center, we have that $c^{(i)}_{km}$ and $d^{(i)}_{km}$ vanish unless $k=1$. In our simulations we choose a coupling strength $\alpha=0.1$, and take the natural frequencies $\omega_i$ to be random with uniform distribution in the interval $[0,2\pi]$ radians per second. Initial conditions are evolved with a fourth order Runge-Kutta integrator with variable step and time series of the phases $\phi_i$ are then collected with a rate of 10 points per second.

To measure the success of the recovery of methods $L_2$  and LASSO,  {we use the measures
\begin{itemize}
\item[] $\#$FP (false positives) consisting of connections that are not present in the true network;
\item[] $\#$FN  \suggestion{(}false negatives\suggestion{)} the connections that were missed by the recovery.
\end{itemize}
We do not take into account the strength of the \suggestion{recovered} connection; instead we simply check whether a certain connection is present or not. We discard connections that are too weak, less than $10\%$ of the largest entry of the coefficient vector}.

\subsubsection{Effects of the length of the time series} 

\begin{figure}[t]
\begin{center}
\includegraphics[scale=0.45,clip]{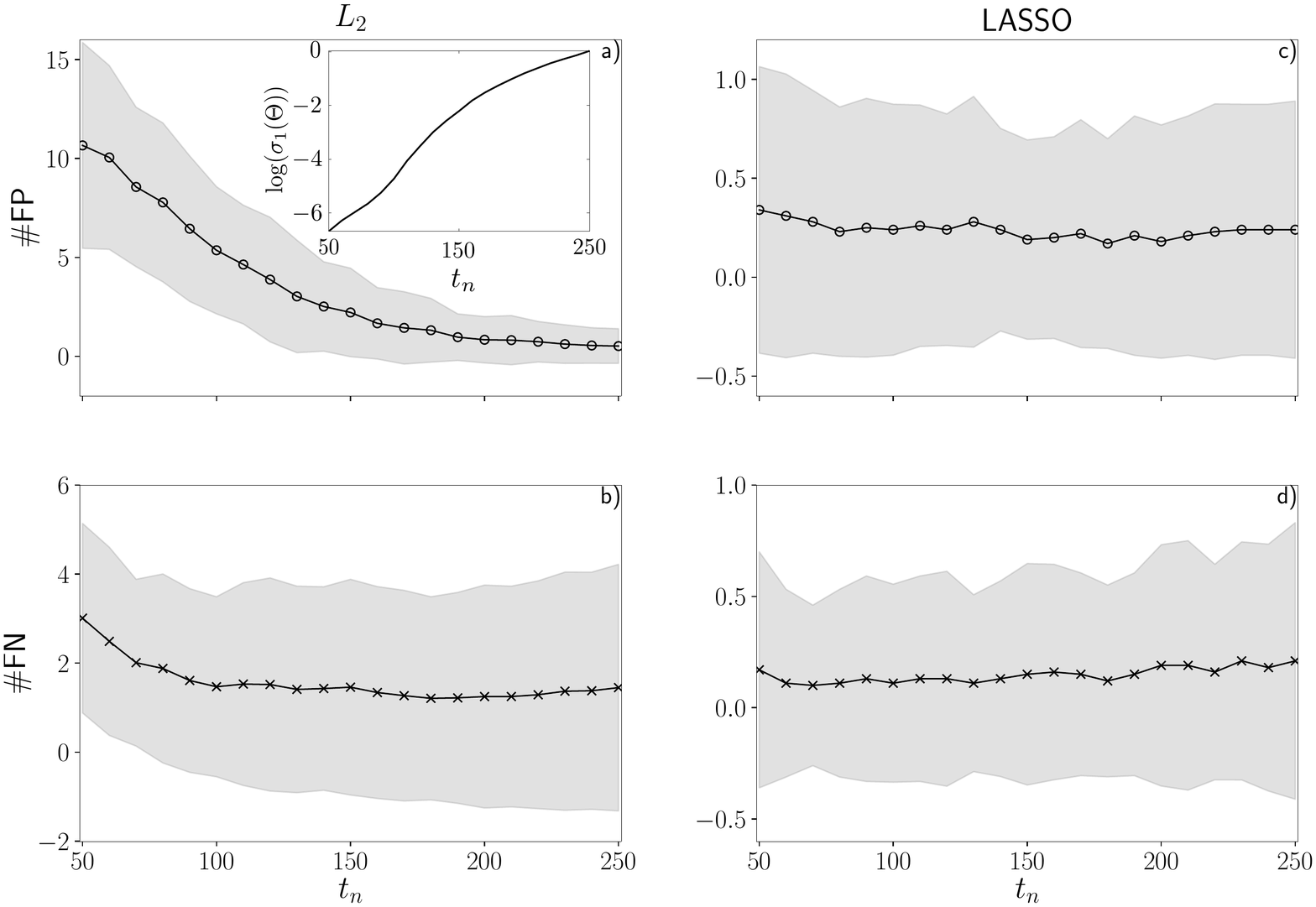}
\end{center}
\caption{  {{\bf Influence of acquisition time $t_n$ on the network recovery}. We consider a directed star graph with $N=10$ and connections diverging from the hub. Panels a) and c) show the false positives $\#$FP (circles) predicted by $L_2$ and LASSO, respectively, as number the acquisition time increases. Panels b) and d) show the false negatives $\#$FN (crosses). Each point is an average over $100$ random initial conditions and the shaded region is the standard deviation. The inset of panel a) shows the logarithm of the minimum singular value of ${\Theta}$, averaged over the $100$ random initial conditions.}} 
\label{fig:CompLASSOL2_Timeincreases}
\end{figure}

\suggestion{In Figure~\ref{fig:CompLASSOL2_Timeincreases},} we show  {$\#FP$ (circles) and $\#FN$ (crosses)}, for the $L_2$ minimization (left column) and for the solution obtained using LASSO (right column), as the acquisition time $t_n$ is varied. These values were averaged over $100$ random initial conditions of our network system with $N=10$ nodes. The LASSO solution is excellent for all values of $t_n$. The $L_2$ minimization performs relatively well if $t_n$ is large, but for small values of $t_n$ it predicts many wrong connections. Similar results were obtained by Napoletani and Sauer~\cite{Domenico_2008}.  

As discussed in the Section~\ref{sec:L2_unstable}, the performance of $L_2$ minimization as a function of $t_n$ seems to be related to $\sigma_1({\Theta})$, the smallest singular value of the matrix ${\Theta}$, which can be small for small $t_n$, as shown in the inset. Subsequently, in Section~\ref{sec:L1_stable_basis_extension}, we show the reason the LASSO approximation is not affected as much by the poor conditioning of ${\Theta}$. \par 

\subsubsection{Effects of the size of the network with fixed length of time series}
\begin{figure}[t]
\begin{center}
\includegraphics[scale=0.4,clip]{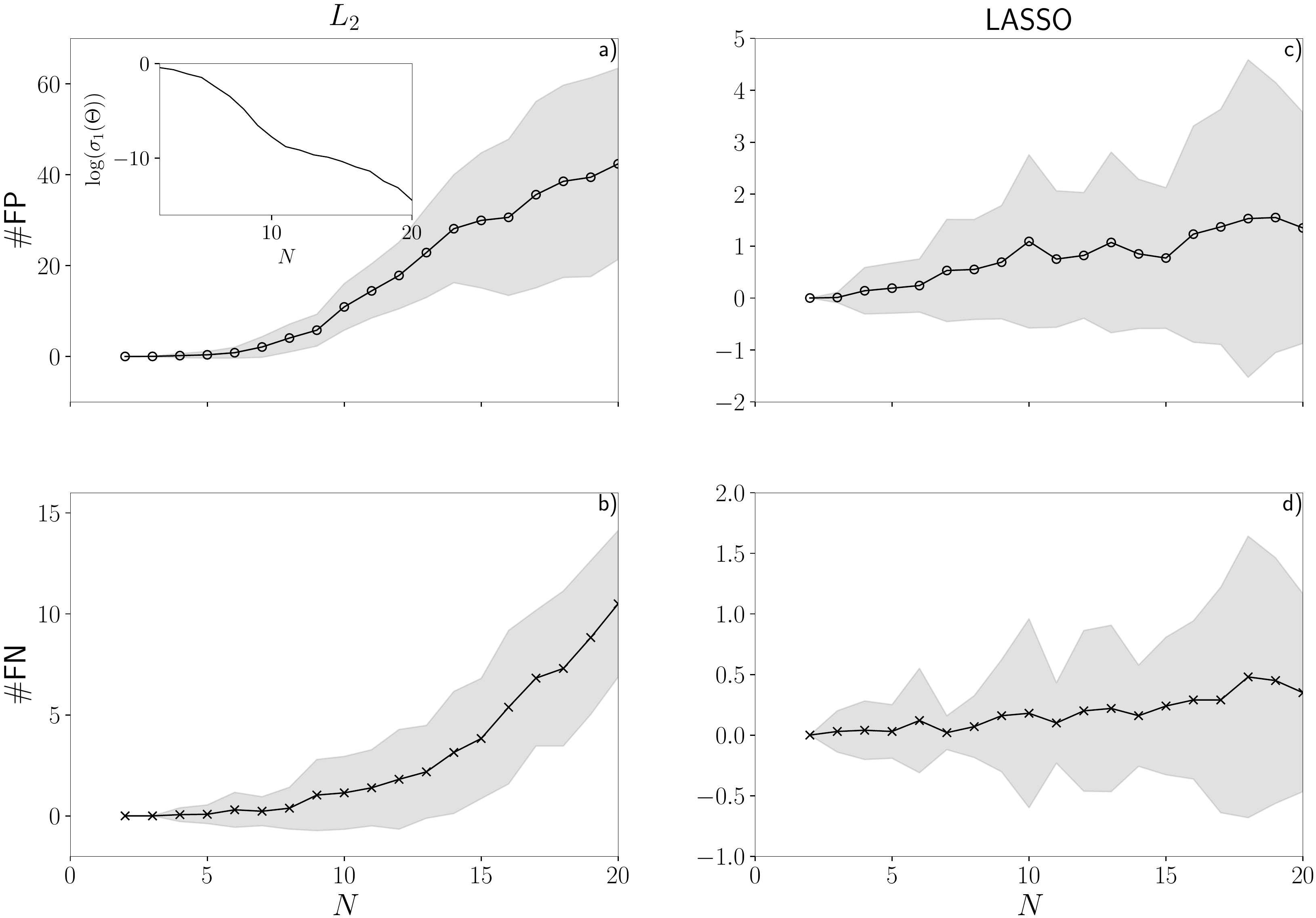}
\end{center}
\caption{  {{\bf Influence network size on the recovery.} We consider a directed star graph of $N$ nodes with connections diverging from. We fix acquisition time $t_n=100$. Panels a) and c) show the false positives $\#$FP (circles) recovered by $L_2$ and LASSO, respectively, as number of nodes $N$ increases. Panels b) and d) show the false negatives $\#$FN (crosses). Each point results from average over $100$ random initial conditions and the shaded region is the standard deviation. In the inset of panel a) we show logarithm of the minimum singular value of ${\Theta}$, averaged over the $100$ random initial conditions.}} 
\label{fig:FixedLengthTimeSeries}
\end{figure}

 {For the directed star graph, in Figure~\ref{fig:FixedLengthTimeSeries} we show $\#FP$ (circles) and $\#FN$  (crosses) as functions of the total number of nodes, for both solutions of $L_2$ minimization and LASSO.} 
The LASSO solution is stable while $L_2$ minimization is accurate for  small networks, with $N\le 7$, and is not able to handle the large-but-sparse configuration. In the inset, we show the corresponding average value of $\log (\sigma_1({\Theta}))$. It suggest a correlation between between the poor performance of the $L_2$ minimization and ill condition of $\Theta$ captured by a small singular value. \par

\subsection{Results for other networks}

\begin{figure}[t]
\begin{center}
\includegraphics[scale=0.45,clip]{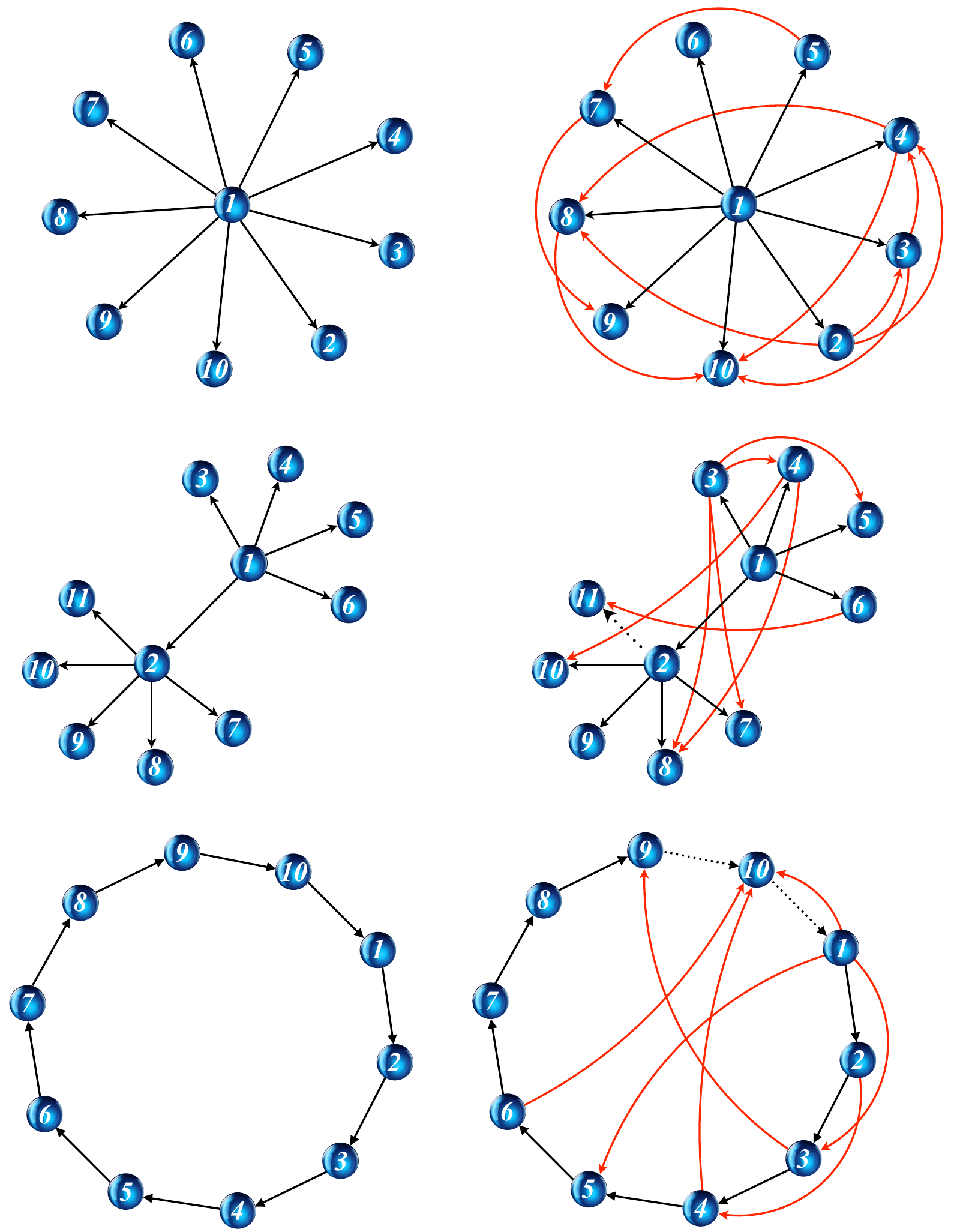}
\end{center}
\caption{{\bf Comparison between LASSO and $L_2$ minimization for three paradigmatic networks}. In the upper panel the directed star, in the mid panel 
one connected directed star forcing another directed star, in the bottom panel a directed ring. The network recovered by the LASSO is presented in the left and shows perfect recovery and the $L_2$ minimization recover is presented in the right. Spurious connections are shown as thin red lines, missing connections are shown as dotted lines. We used a  single random initial condition in each case and  $t_n=100$. } 
\label{fig:CompLASSOL2}
\end{figure}

In this section we briefly consider some other sparse networks: \suggestion{the twin stars, which consists of two stars joined by a single link, and a ring,} both are illustrated in Figure~\ref{fig:CompLASSOL2}. \incorrect{In the twin stars configuration node $1$ drives nodes $2$ to $6$, while node $2$ drives nodes $7$ to $11$. In the ring configuration node $i$ drives its following neighbour $i+1$, and node $N$ drives node $1$.}

We again use $\alpha=0.1$ and take the natural frequencies $\omega_i$ to be random with uniform distribution in the interval $[0,2\pi]$. Initial conditions are evolved with a Runge-Kutta integrator and time series of the phases $\phi_i$ are then collected with time steps of $0.1$.

In Figure 4 we show the connections that were recovered by the two methods, $L_2$ minimization and LASSO, from a single random initial condition propagated for $t_n=100$. We performed a kind of hard thresholding, by discarding connections that were too weak (we considered coupling strengths smaller that 10\% of the largest one to be weak). 

The results from LASSO are excellent in all cases, but the $L_2$ minimization does not perform so well: it  fails both to recover existing connections (false negatives depicted with dotted lines in Figure 3) and recovers false positive (thin grey lines).

\subsection{Effects of basis extension}

 {We discuss how the inclusion of new functions in the basis can affect the recovery.  Our first example is shown in Figure \ref{fig:FixedLengthTimeSeries}. Since the network is a directed star with connections diverging from the hub, the recovery of each node is independent as the hub acts as a master to the leaves. Thus, increasing the network size and recovering the connections of a given node has the same effect as including new ({\it a posteriori}) unnecessary functions in the basis. We could wrongly expect this basis extension would not influence the recovery method. Figure \ref{fig:FixedLengthTimeSeries} shows that the $L_2$ recovery is strongly affected by such extensions as the inclusions of new functions, while keeping the length of the time-series fixed,  makes the operator $\Theta$ ill-conditioned.} 

\begin{figure}[!th]
\begin{center}
\includegraphics[scale=0.3,clip]{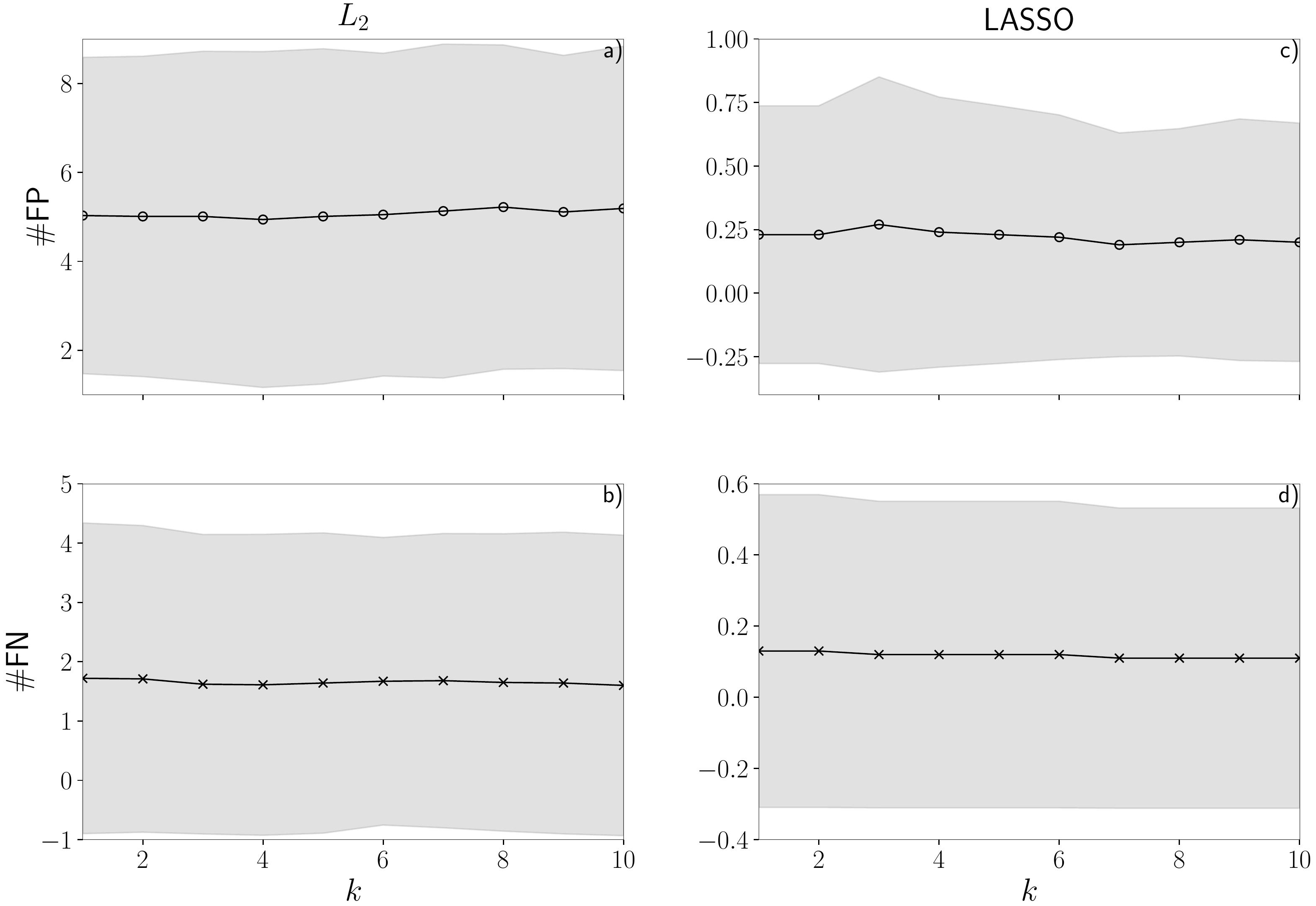}
\end{center}
\caption{  {{\bf Effects of basis extension in the recovery}. 
We consider a fixed acquisition time $t_n = 100$ and a directed ring with $N=10$.
First, we recover the network without higher harmonics in the phases corresponding to $k$. Then, we extend the basis to include higher harmonics $\mathcal{E}_k = \{\sin (m \vartheta_k), \cos(m \vartheta_k)\}_{m = 2}^{10}$ of a node $k$ iteratively. Thus, for each $k$ we include $16$ new functions in the basis and apply the recovery methods while keeping the previously added functions. Panels  a) and c) show the false positives $\#$FP (circles) recovered by $L_2$ and LASSO, respectively, as a function of $k$. Panels b) and d) show the false negatives $\#$FN (crosses). 
Each point is an average over $100$ initial conditions and the shaded region is the standard deviation.}} 
\label{fig:CompLASSOL2_grow_basis}
\end{figure}

{   Next, we notice that the function $g$ in Eq.(\ref{eq:phase_model}) plays no role in the dynamics when phases are slow variables. We study the effect of the inclusion of such functions in the recovery process. Figure~\ref{fig:CompLASSOL2_grow_basis} shows the results of such basis extension for a directed ring. The basis extension is made using higher harmonics $\mathcal{E}_k = \{\sin (m \vartheta_k), \cos(m \vartheta_k)\}_{m = 2}^{10}$ for each node $k$. 
Starting from $k=0$ we include the new functions of a node $k$ while keeping all previously included functions. Thus, in the first iterator $k=1$ we include $16$ new functions and at the end of the process $k=10$ we include $160$ functions.   
We observe that the number of false positives and negatives remains unaltered as the basis is increased either for $L_2$ and LASSO. We notice that LASSO remains stable under basis extension.}

\section{Stability of sparse networks under basis extension}

\subsection{$L_2$ is unstable under basis extension}
\label{sec:L2_unstable}

When we extend the basis, probing new possible connections, we face a problem as $[A,B]$ may have small singular values, leading to instabilities. This means that even if 
\begin{linenomath*}
\begin{equation}
A x = b \nonumber
\end{equation}
\end{linenomath*}
\noindent
has a sparse solution, it may happen that 
\begin{linenomath*}
\begin{equation}
[A,B]w = b+z \nonumber
\end{equation}
\end{linenomath*}
\noindent
has a solution $w$ that is far from being sparse in its restriction to the components corresponding to $x$, here \suggestion{$z$ captures small measurement errors.}
This would mean that the basis extension is unstable. 
 
{\suggestion{Our next proposition characterizes this situation}. We prove it using the concept of principal angle between subspaces, in particular the largest principle angle between the orthogonal complement of the image of matrix $A$, $(\im~A)^{\perp}$, and the image of the matrix $B$, $\im~B$.

\begin{proposition}\label{prop:L_2_instability}
Let ${A}\in \mathbb{R}^{n \times p}$ be a column full rank matrix,  $b \in \im~A$ and $z \in (\im~A)^{\perp}\backslash\{0\}$. Let $x^*$ be the unique solution of the problem 
\begin{linenomath*}
\be
\min_{x\in \mathbb{R}^{p}} \|  A x - b - z \|_2. \nonumber
\ee
\end{linenomath*}
Let ${B} \in \mathbb{R}^{n \times q}$ be such that the matrix \suggestion{concatenation}
$C = [A,B]$ is also column full rank with $n > p + q$. Let $r = \min\{p, q\}$ and the principal angles between the subspaces $(\im~A)^{\perp}$ and $\im ~B$ satisfy: $0 < \beta_1 < \dots < \beta_r < \frac{\pi}{2}$. Let $\hat w = (\hat w_1, \hat w_2)$ be the unique solution of the problem 
\begin{linenomath*}
\be
\label{eq:partitioned_l2_minization}
\min_{w \in \mathbb{R}^{p+q}} \| C w - b - z \|_2.    \nonumber
\ee
\end{linenomath*}
Then for a generic $z>0$ given a natural number $N_0 > 0$ there is a $\varepsilon>0$ such that if 
$| \beta_r - \pi/2| < \varepsilon$ we obtain 
$\|  x^* - \hat w_1 \|_2 > N_0$. 
\end{proposition}

We prove this proposition in Appendix \ref{Insta}.
}
The above proposition explains the instability we observed in the numerical results, which are also in agreement with the observations made by Napoletani and Sauer~\cite{Domenico_2008}. \par 

As a remark, when the dynamics is chaotic the columns of the matrix $\Theta$ behave as pseudorandom vectors. Let us assume that $p=q$ for the matrices $A$ and $B$. 
Thus we can think of the column spaces of $A$ and $B$ as two $p$-dimensional vector spaces taken at random from a larger $n$-dimensional space, $n > 2p$. The principal angles between them have a joint multivariate beta distribution \cite{absil2006largest}; from well known random matrix theory results, it then follows that, as $n \rightarrow \infty$ with $p = \xi n$, the average value of the smallest principal angle satisfies $\cos(\beta_1) = 4 \xi(1-\xi)$. The value $\xi \rightarrow 1/2$ corresponds to the case $[A,B] = \Theta$, when the principal angles tends to $0$. This indicates that, in the large basis limit, instability is generic.  

\subsection{Sparse Solutions}

To fix the problem of basis instability, we take into account the sparsity of the network, i.e. the fact that only a few of the coefficients we are looking for will be nonzero. 
Empirically, when we take this into account we can drastically reduce the number of data points needed 
for the reconstruction as well as gain the stability of the reconstruction starting from a seed. 
First, we have some definitions.

\begin{definition} A vector is said to be $s$-sparse if it has at most $s$ nonzero entries
\begin{linenomath*}
\begin{equation}
\| x \|_{0} \le s.
\end{equation}
\end{linenomath*}
\end{definition}

\noindent
\suggestion{{\bf Notation.}  $x_s$ is the vector obtained from $x$ when all but the $s$ largest entries are set to zero.}

\begin{definition}For each positive integer $s$, define the $s$th restricted isometry constant $\delta_s = \delta_s(A)$ of a matrix $A$ 
as the smallest number such that
\begin{linenomath*}
\begin{equation}\label{RIP_p}
(1 - \delta_s) \| x\|_2^2 \le \| A x\|_2^2 \le (1 + \delta_s) \| x\|_2^2
\end{equation}
\end{linenomath*}
for all $s-$sparse vectors $x$. \suggestion{
The matrix $A$ which has $\delta_s \in (0, 1)$ is said to satisfy the restricted isometry property (RIP).}
\end{definition}

Assume that we have measurements corrupted with noise so that
\begin{linenomath*}
\begin{equation}
Ax = b + z,
\end{equation}
\end{linenomath*}
where $z$ is an unknown noise term. In this context, \suggestion{one may reconstruct} $x$ as the solution to the convex optimization problem
\begin{linenomath*}
\be \label{eq:quadractically_constrained_basis_pursuit}
\min_{\tilde x \in \mathbb{R}^p} \| \tilde x\|_{1} \mbox{~subject to ~} \| A\tilde x - b \|_2 \le \varepsilon
\ee
\end{linenomath*}
where $\varepsilon$ is an upper bound on the noise. 
The next statement shows that one can stably reconstruct $x$ as long as the matrix $A$ has a controlled restricted isometry constant $\delta_s$.
\begin{theorem}[Noisy recovery - Theorem 1.3 \cite{CANDES_2008}]\label{candes}
Assume that $\delta_{2s} < \sqrt{2}-1$ and $\| z \|_2 \le \varepsilon$. 
Then the solution $x^*$ to Eq. (\ref{eq:quadractically_constrained_basis_pursuit})
 satisfies
 \begin{linenomath*}
\begin{equation}
\|x^* - x\|_2 \le C_0 s^{-1/2} \| x  - x_s\|_1 + C_1 \varepsilon
\end{equation}
\end{linenomath*}
for some constants $C_0$ and $C_1$.
\end{theorem}

The proof of this result can be found in \cite{candes1,CANDES_2008}. 
Thus the major issue is whether we can find a set of basis functions which yields good properties such as RIP for the matrix $\Theta$. We suggest that the set of basis functions must be built over dynamical information from the underlying dynamical system generating the time series. A key property here is the \emph{coherence} of a matrix
\begin{definition}[Coherence]
Let $A \in \mathbb{R}^{n \times m}$ be a matrix with $L_2$-normalized columns $v_1,\dots,v_m$ 
its coherence $\eta(A)$ is defined as
\begin{linenomath*}
\be
\eta = \max_{i\not = j} | \langle v_i, v_j \rangle |. \nonumber
\ee
\end{linenomath*}
\end{definition}
The coherence quantifies the linear independence of pairs of matrix columns. Consequently, it is intrinsically linked to RIP constant $\delta_s$. This will play essential role in Section~\ref{sec:L1_stable_basis_extension} to prove the stability of the $L_1$ minimization problem under basis extension. 

\subsection{Basis Adaptation guarantees coherence}
Let $\mathbb{T} = \mathbb{R}/2\pi\mathbb{Z}$ be the torus. From here on our theoretical formulation and analysis is described in terms of a map denoting the dynamics. This assumption is not harmful since the phase dynamics recovery on $\mathbb{T}^N$ is given by the time-one map $f$ of the flow. This map is induced by the Euler approximation of the differential equations and the sampling procedure of the trajectories.  \par 
In the following exposition we will denote by $X$ the metric space being either 
a compact subset of $\mathbb{R}^d$ or a parallelizable manifold such as the 
torus $\mathbb{T}^d$. We assume the map denoting the dynamical system is $C^{r}(X)$ with $r \geq 1$. This will contain all examples in the paper and avoid a technical detour.  {We denote $\psi$ as basis functions representing the map $f$ and the functions $\varphi$ and $\phi$ are observables.} We understand sparse representation of the map as
\begin{definition}[Sparse Representation]
Let $f : X \rightarrow X$ and $\mathcal{L} = \{ \psi_i \}_{i=1}^m$ be a set of basis functions 
with $\psi_i: X \rightarrow X$ for $i = \{ 1,\dots, m\}$ such that
\begin{linenomath*}
\be
f = \sum_{i=1}^{m} c_i \psi_i. \nonumber
\ee
\end{linenomath*}
We say that $f$ has an $s$-sparse representation in $\mathcal{L}$ if the vector $x = (c_1, \dots, c_m)^{\dagger}$ is $s$-sparse. 
\end{definition}
\textbf{Dynamical information: Ergodicity.} We focus our analysis on ergodic dynamical systems. A well-known property is that the time average of an observable evaluated at a typical orbit converges to the space average. This is more generally stated in the following Theorem~\ref{thm:Birkhoff_theorem}

\begin{theorem}[Birkhoff Ergodic theorem]\label{thm:Birkhoff_theorem}
 Let $(f, \mu)$ be a discrete ergodic dynamical system on the compact metric space $X$. Given any $\varphi \in L_1(\mu)$, there exists a set of initial conditions $E \subset X$ with $\mu(E) = 1$ such that for any $\varepsilon > 0$ and $x_0 \in E$ there exists $n_0(\varepsilon, x_0) > 0$ where the following holds
 \begin{linenomath*}
 \begin{equation}\label{eq:time_and_space_average}
     \Big|\frac{1}{n} \sum_{k = 0}^{n - 1} \varphi \circ f^k(x_0) - \int \varphi d\mu\Big| < \varepsilon \qquad \forall~ n > n_0.
 \end{equation}
 \end{linenomath*}
 \end{theorem}
Birkhoff Ergodic theorem is the main ingredient to calculate the coherence for ergodic dynamical systems in Theorem~\ref{thm:coherence}. It introduces a change of inner product: instead of looking at the euclidean inner product among vectors on $\mathbb{R}^n$, we approximate it by the inner product on the space of integrable functions with respect to the ergodic measure. \par 
Theorem~\ref{thm:coherence} goes beyond. It states that for any discrete ergodic dynamical system whose measure has density, we can construct a set of basis functions adapted to the ergodic measure via Gram-Schmidt procedure\footnote{These adapted basis functions are related to the Bounded Orthonormal System (BOS) in Foucart and Rauhut~\cite{Foucart_mathematical_compressive_sensing}. They differ in respect to the choice of the reference measure. BOS carries the measure given by the uniform sampling procedure whereas here it comes along the observed trajectory.}. These adapted basis functions do not harm the sparsity representation of the map and has control over the coherence of the matrix for large enough data. \par 

Our result is related to what was obtained by Tran and Ward~\cite{Tran_Ward_CLT_Lorenz}. The authors use Central Limit Theorem applied for Lorenz systems perturbed over time to obtain the null-space property (which is a weaker property than RIP) for a similar version of the matrix $\Theta$. \par 

To our best knowledge, our results are one of the few examples to advance the search for basis functions adapted to the dynamical system generating the time series. Recently, Hamzi and Owhadi~\cite{hamzi2020learning} proposed a kernal-based method in a similar direction. \par 

\emph{Drawback.} It is worth noting that Theorem~\ref{thm:coherence} is an existence statement since it requires that the sparse representation of the dynamical system is known a priori. Besides it is valid for large enough data. To determine the minimum amount of data for controlling the coherence of $\Theta$, it would be necessary to know the speed of convergence of the Birkhoff sums for the basis functions. This will be done in the near future.

 \begin{theorem}[Ergodic Coherence]\label{thm:coherence}
Let $(f, \mu)$ be an ergodic dynamical system with $\mu$ absolutely continuous with respect to Lebesgue (${\rm Leb}(X)$). Let $\mathcal{L}_0$ be a set of basis functions such that $f$ has an $s$-sparse representation in $\mathcal{L}_0$. Given $\eta_0 > 0$ and $\varepsilon \in (0, 1)$ there is a set of basis functions $\mathcal{L}$, $n_0>0$ and a good set of initial conditions $G \subset X$ such that 
\begin{itemize}
\item[(i)] $\mu(G) > 1 - \varepsilon$, and for any $x_0 \in G$ and  $n > n_0$  we have $\eta(\Theta(\mathcal{L})) < \eta_0$.
\item[(ii)] the representation of $f$ in $\mathcal{L}$ is also $s$-sparse.
\end{itemize}
 \end{theorem}

 \begin{proof}
We develop the argument assuming that $X \subset \mathbb{R}$. To generalize for large dimensions or for $\mathbb{T}^{d}$ it is enough to break down the problem in terms of coordinates.  
The main ingredient in the proof is the  Birkhoff's Ergodic Theorem \ref{thm:Birkhoff_theorem}. 
Having the ergodic theorem we split the proof in three steps. 

{\it Step 1: Ergodicity and basis adaptation.} Let $\mathcal{L}_0 = \{ \psi_1, \dots, \psi_m \}$ \suggestion{ be a set of basis functions, where each $\psi_i: X \to X$}. We perform a Gram-Schmidt process in $L_2(\mu)$ and obtain an orthogonal basis 
\begin{linenomath*}
\be    
\hat{\mathcal{L}} = \{ \varphi_1, \dots, \varphi_m\}. \nonumber
\ee
\end{linenomath*}
Notice that since $\mu = \nu {\rm Leb} $ we define 
\begin{linenomath*}
\be
\phi_i = a_i \varphi_i  \nonumber
\ee
\end{linenomath*}
where   {$a_i^2 = 1/ \int \varphi_i^2 ~\nu ~d {\rm Leb}$} such that  $\mathcal{L} = \{ \phi_i \}_{i=1}^{m}$ is an orthonormal system with respect to   {$L_2(\mu)$} in the span of $\mathcal{L}_0$.
For an arbitrary initial condition $x_0$, let
 \begin{linenomath*}
    \be
    u_i :=  
       \frac{1}{\sqrt{n}} \begin{pmatrix}
           \phi_i(x_0) \\
           \vdots \\
           \phi_i(f^{n - 1}(x_0)) \\
        \end{pmatrix} \qquad u_j :=
       \frac{1}{\sqrt{n}} \begin{pmatrix}
           \phi_j(x_0) \\
           \vdots \\
           \phi_j(f^{n - 1}(x_0)) \\
        \end{pmatrix} 
    \ee
    \end{linenomath*}
be the $i$th and $j$th columns of the matrix $\Theta(\mathcal{L}) \in \mathbb{R}^{n \times m}$. Then notice that the inner product between columns $i$ and $j$ is
\begin{linenomath*} 
\begin{eqnarray*}
    \langle u_i, u_j \rangle &=& \frac{1}{n} \sum_{k=1}^{n} \phi_i(f^k(x_0)) \phi_j(f^k(x_0))  \\
&=& \frac{1}{n} \sum_{k=1}^{n} (\phi_i \cdot  \phi_j )\circ (f^k(x_0)) \\
&=:& \frac{1}{n} S_n(\phi_i \cdot  \phi_j )(x_0). 
 \end{eqnarray*} 
 \end{linenomath*}
From the smoothness of the map $f$, $(\phi_i \cdot  \phi_j)$ is integrable $L_1(\mu)$ so by Birkhoff Ergodic theorem there is a set $G_{ij}$ such that $\mu(G_{ij})$ has full measure and for each $x_0 \in G_{ij}$ and $\varepsilon_1 > 0$ there is $n_0 > 0$ such that for any $n>n_0$ we have
\begin{linenomath*} 
\begin{eqnarray}
\left| \langle u_i, u_j \rangle -  \int \phi_i \cdot  \phi_j d\mu \right| \le \varepsilon_1  \nonumber \\
\left| \langle u_i, u_j \rangle -  \delta_{ij}\right| \le \varepsilon_1 \nonumber
\end{eqnarray} 
 \end{linenomath*}
where $\delta_{ij}$ is the Kronecker delta. 

{\it Step 2: Large measure of initial conditions for the basis}. 
Hence, we are interested in the subset with cardinality $K = \frac{m(m - 1)}{2}$
\begin{linenomath*}
\be\label{eq:set_of_observables_}
\mathcal{G} = \{ (\phi_{i}\cdot \phi_j) ~|~i, j = 1, \dots, m\} \subset L_{1}(\mu)
\ee
\end{linenomath*}
{where each element corresponds to pairwise multiplication of basis functions in $\mathcal{L}$. We aim at finding a good set $G$ of initial conditions where the control of $n_0$ is uniform.}
    
{Using Egoroff's theorem~\cite{folland2013real} we can make the Birkhoff sum $\frac{1}{n}S_n \phi$ converge uniformly on a large measure set $G_{\phi}$ of $X$ instead of the ``almost every point" convergence. }
{Fix $\eta_0 > 0$ and take $\varepsilon/(2K)$. For each observable $\phi$ in the set of Equation (\ref{eq:set_of_observables_}), the precision $\varepsilon/(2K)$ determines a subset $G_{\phi}$ of $X$ which by Egoroff's theorem has measure $\mu(G_{\phi}) > 1 - \frac{\varepsilon}{2K}$ where the convergence of $\frac{1}{n}S_n \phi$ is uniform. So, we take the set of initial conditions as 
 \begin{linenomath*}
\be G = \bigcap_{\phi \in \mathcal{G}} G_{\phi}.
\ee
\end{linenomath*}
}Using the complement of $G$ we can calculate that 
 \begin{linenomath*}
\be
    \mu(G)  > 1 - \varepsilon. \nonumber
 \ee
 \end{linenomath*}
This determines the set of initial conditions for which we can calculate the coherence of the matrix $\Theta(\mathcal{L})$. Due to uniformity of initial conditions in $G$, for each observable $(\phi_{i}\cdot \phi_j)$ in the set of Equation~\ref{eq:set_of_observables_} there exists $n_{i, j} > 0$ such that the inner product of any two distinct normalized column vectors $|\langle v_i, v_j\rangle|$ has the following form for any $n > n_{i, j}$
\begin{linenomath*}
\be
        |\langle v_i, v_j\rangle| = \Bigl|\frac{1}{n} S_n (\phi_i \cdot \phi_j)  (x_0)\Bigr| \leq \eta_0. \nonumber
 \ee
  \end{linenomath*}
Take $n_0 := \max_{i \neq j} n_{i, j}$ and this proves the statement.   

 {\it Step 3: Sparsity.}   
Thus we are only left to prove sparsity. We know by assumption that there is a sparse 
solution to 
 \begin{linenomath*}
\be
\Theta(\mathcal{L}_0) x_s = v.\nonumber
\ee
 \end{linenomath*}

Let us rearrange $\mathcal{L}_0$ such that $x_s$ has only its first $s$ entries nonzero.
Next, the Gram-Schmit process reduces to a QR decomposition that is 
 \begin{linenomath*}
\be
\Theta(\mathcal{L}_0) = \Theta(\mathcal{L}) R     \nonumber
\ee
 \end{linenomath*}
\noindent
thus, 
 \begin{linenomath*}
\be
\Theta(\mathcal{L}_0) x_s = \Theta(\mathcal{L}) \hat x_s  \nonumber
\ee
 \end{linenomath*}
\noindent
where 
\begin{linenomath*}
\be
\hat x_s = R x_s \nonumber
\ee
\end{linenomath*}
but $R$ is upper triangular and thus by construction only the first $s$ entries of $\hat x_s$ will be nonzero.
\end{proof}

\subsection{Sparse Solutions are stable under basis extension}
\label{sec:L1_stable_basis_extension}
Next, we wish to prove that once the basis is adapted and the initial expert
guess is meaningful, extending the basis is not harmful for the solution. 
First, we need the following result relating the coherence and restricted isometry constant of a matrix

\begin{proposition}\label{CoRIP}
If the matrix $A$ has $L_2$-normalized columns, then its RIP constant satisfies
\begin{linenomath*}
\be
\delta_s \le \eta (s-1), \,\,\,s\ge2 \nonumber
\ee
\end{linenomath*}
\end{proposition}
\begin{proof}
See Foucart and Rauhut~\cite[Prop. 6.2, p.134]{Foucart_mathematical_compressive_sensing}.
\end{proof}

Next proposition proves that, given a  set of basis functions which represents $f$ sparsely, the minimization problem from Cand\`es Theorem~\ref{candes} has a solution that approximates the true sparse solution. Moreover, using Theorem~\ref{thm:coherence}, which introduces a orthonormal set of basis functions $\mathcal{L}$ and a matrix $\Theta(\mathcal{L})$, we can find a sub-matrix of $\Theta(\mathcal{L})$, $A(\mathcal{L})$, which approximates the same solution in a smaller space. \par 
It is worth noting that both LASSO and quadratically constrained basis pursuit are $L_1$ minimization problems related to each other.  More precisely, for each solution $x^{\star}$ of LASSO there exists a $\varepsilon := \varepsilon_{x^{\star}} > 0$ such that $x^{\star}$ is solution of Equation~(\ref{eq:quadractically_constrained_basis_pursuit}), see Fourcart and Rauhut~\cite[Proposition 3.2]{Foucart_mathematical_compressive_sensing}. So, our results using the quadratically constrained basis pursuit  are extended to LASSO solutions as well.\par 

\begin{proposition}[Sparsity level is attained]\label{prop:A_L_exists}
Let $\mathcal{L}_0$ be a set of basis functions with cardinality $m$ such that $f$ has a $s$-sparse representation in $\mathcal{L}_0$. Then, there is $n_0 > 0$, a large set of initial conditions and a basis $\mathcal{L}$ such that we find a matrix $A(\mathcal{L}) \in \mathbb{R}^{n \times p}$ where $s < p < m$ and the solution $x^*$ of the reconstruction problem 
\begin{linenomath*}
\be \label{eq:A_L_quadractically_basis_pursuit}
\min_{\tilde{x} \in\mathbb{R}^p} \| \tilde x\|_1 \mbox{~subject to~} \| A(\mathcal{L}) \tilde x - v\|_2 <\varepsilon
\ee
\end{linenomath*}
attains the sparse representation of $f$.
\end{proposition}
\begin{proof}
Using Proposition~\ref{CoRIP} together with Theorem~\ref{thm:coherence} we  conclude the following: let $1 < s < m$ be the sparsity level of the representation of the map $f$ with respect to the proposed set $\mathcal{L}_0$.  \par 
By assumption we know there exists a sparse solution $x_s \in \mathbb{R}^{m}$ such that $\Theta(\mathcal{L}_0)x_s = v$. We rearrange $\mathcal{L}_0$ such that $x_s$ has only its first $s$ entries nonzero. 
Fix $0 < \eta_0 < (\sqrt{2}-1)/
(2s - 1)$. By Theorem~\ref{thm:coherence} there exists an orthogonal basis $\mathcal{L}$, $n_0>0$ and a large set of initial conditions that 
$
\eta(\Theta(\mathcal{L}))\le \eta_0
$
and $\hat{x}_s \in \mathbb{R}^{m}$. Thus from Proposition \ref{CoRIP} 
\begin{linenomath*}
\be
\delta_{2s}(\Theta(\mathcal{L})) < \sqrt{2} - 1.
\ee
\end{linenomath*}
By Theorem~\ref{candes}, the sparse solution $\hat{x}_s$ is approximated by the solution of the quadractically constrained basis pursuit problem. \par 
Let $p, q \in \mathbb{N}$ be chosen such that $s < p, q < m$ and $m = p + q$. Without loss of generality, we can rearrange the basis elements in such way $\Theta(\mathcal{L}) = [A(\mathcal{L}), B(\mathcal{L})]$ where 
$A(\mathcal{L}) \in \mathbb{R}^{n \times p}$ and $B(\mathcal{L}) \in \mathbb{R}^{n \times q}$. Moreover, using  $A(\mathcal{L})$ in the quadratically constrained basis pursuit problem the solution approximates the sparse solution $\hat{x}_s$ through a vector lying in $\mathbb{R}^{p}$ . This is true because $\delta_{2s}(\Theta(\mathcal{L}))$ is an upper bound for $\delta_{2s}(A(\mathcal{L}))$ and $\delta_{2s}(B(\mathcal{L}))$. \par 
For the noiseless case we could say that $A(\mathcal{L})$ is the minimum matrix such that the minimization problem attains the sparse solution. \par 
\end{proof}

The existence of a sub-matrix of $\Theta(\mathcal{L})$ in the above proposition indicates that we can use Theorem~\ref{thm:coherence} in a different way to guarantee that sparse solutions are stable under basis extension. The following corollary states this stability more precisely. \par 

\begin{coro}[Stability under basis extension]\label{coro:L1_stable}
 Suppose $\mathcal{L}_0$ is a subset of basis functions with cardinality $p < m$ such that $f$ has a $s$-sparse representation in $\mathcal{L}_0$. Denote $x_s \in \mathbb{R}^{p}$ the unique sparse solution of Equation (\ref{eq:A_L_quadractically_basis_pursuit}). 
 Then there is $n_0 > 0$, a large set of initial conditions and a basis $\mathcal{L}$ 
such that $w = (x_s, 0) \in \mathbb{R}^{m}$ is solution of
\begin{linenomath*} \be
 \Theta(\mathcal{L}) w = [A(\mathcal{L}), B(\mathcal{L})] w = v \nonumber
\ee
\end{linenomath*}  
and the solution $w^{*} = (w_1^{*}, w_2^{*})$ of
\begin{linenomath*} \be
    \min_{\tilde{w} \in \mathbb{R}^{m}} \|\Tilde{w}\|_1 \mbox{~subject to~} \|[A(\mathcal{L}), B(\mathcal{L})] \tilde{w} - v\|_2 < \varepsilon \nonumber
 \ee \end{linenomath*}
satisfies
\[ \label{eq:split_estimate_Candes}
\|w_1^* - x_s\|_2 \le C \varepsilon \quad \mbox{and} \quad  \|w_2^*\|_2 \le \hat{C} \varepsilon 
\]
for constants $C$ and $\hat{C}$.
\end{coro}
\begin{proof}
Thinking in the reverse direction as in the previous proposition we could assume $\mathcal{L}_0$ is a subset of basis functions with cardinality $p < m$ such that $f$ has a $s$-sparse representation in $\mathcal{L}_0$. Then by Theorem~\ref{thm:coherence}, Proposition~\ref{CoRIP} and \incorrect{Candès} theorem~\ref{candes} there is $n_1 > 0$, a large set of initial conditions and a basis $\mathcal{L}$ such that $A(\mathcal{L}) \in \mathbb{R}^{n \times p}$ satisfies Equation (\ref{eq:A_L_quadractically_basis_pursuit}). \par 
The key fact is the finiteness of the set of basis functions. Let us denote by $\mathcal{L}_0^{c}$ the complement of $\mathcal{L}_0$. If we take the union $\mathcal{L} \cup \mathcal{L}_0^{c}$ we can apply Theorem~\ref{thm:coherence} for this set.  
Since $\mathcal{L}$ is already orthonormal, the Gram-Schmidt procedure is necessary only for the functions of $\mathcal{L}_0^{c}$. Adjusting $n_0 > n_1 > 0$ and the initial conditions, and using orthonormality we can guarantee continuity of the unique sparse solution of Equation (\ref{eq:A_L_quadractically_basis_pursuit}) in the larger space.  
The estimate in Equation~(\ref{eq:split_estimate_Candes}) is given by applying Theorem~\ref{candes}.  
\end{proof}

\section{Conclusions}

We  considered the problem of recovering, from phase dynamics, the interaction structure of a sparse network of oscillators. We  compared two different recovery methods, both based on a Fourier expansion of the interaction functions. One of them is the traditional least squares approximation, which finds the vector of coefficients that minimize the $L_2$ error of the approximation and has been successful in previous approaches. The other is LASSO.  
For small networks and when long data sets are available, both approaches are equivalent. But we have found LASSO to be much more apt to sparse network configurations and short times 
than the $L_2$ minimization.  We showed that LASSO can perform remarkably well when dynamical information is taken into account and the basis functions \incorrect{are} adapted. This adaptation leads to unique solutions to the minimization problem that are also stable under basis extension.  Once the basis is adapted to the dynamics LASSO recovers sparse networks with excellent precision even when only relatively little data is available.  \\

\noindent
{\bf Acknowledgments}  M.N. was supported as Visiting Professor by FAPESP grant 018/22503-8. E. R. S. was supported by FAPESP grant 2018/10349-4. T.P. was supported by FAPESP Cemeai grant 2013/07375-0, by Royal Society grant NAF$\backslash$R1$\backslash$180236 and  by the Serrapilheira Institute grant Serra-1709-16124. We thank Sebastian van Strien, Jeroen Lamb, Dmitry Turaev, Richard Cubas and Zheng Bian for useful discussions.

\appendix 

\section{Stability of Lasso under noise}

We analyse the effect of noise by adding, to the original equations of motion, Eq. (1), a term $\sqrt{2D} \dot B_i (t)$, with a homogeneous complex Wiener process  $B_i(t) = \xi_i (t)  + \rm j \zeta_i(t)$ with $\langle \xi_k(t)\xi_i(s) \rangle = \langle \zeta_k(t)\zeta_i(s) \rangle = \delta_{ik} \delta(t-s)$, and $D=$diag$(\eta,\eta)$ and we obtain
\be \label{phasenoise}
\dot{\phi_i}=\omega_i+\alpha\sum_{k=1}^NC_{ik} \sin(\phi_k-\phi_i) + \mathcal{N}_i(\phi_i,t)
\ee
where the noise term is given by 
\be \label{eq:Wiener_process}
\mathcal{N}_i(\phi,t) = \sqrt{2 \eta}(\cos(\phi)\xi_i(t)+\sin(\phi)\zeta_i(t)).
\ee
The noisy equations of motion are integrated by Euler's method with a time step of 0.1, the time series for the phases are obtained by Hilbert transform and a Savitzky-Golay filter is applied to them, before the time derivative is calculated. The filtered phases are then used in the matrix $\Theta$. \par
\begin{figure}[t]
\begin{center}
\includegraphics[scale=0.57,clip]{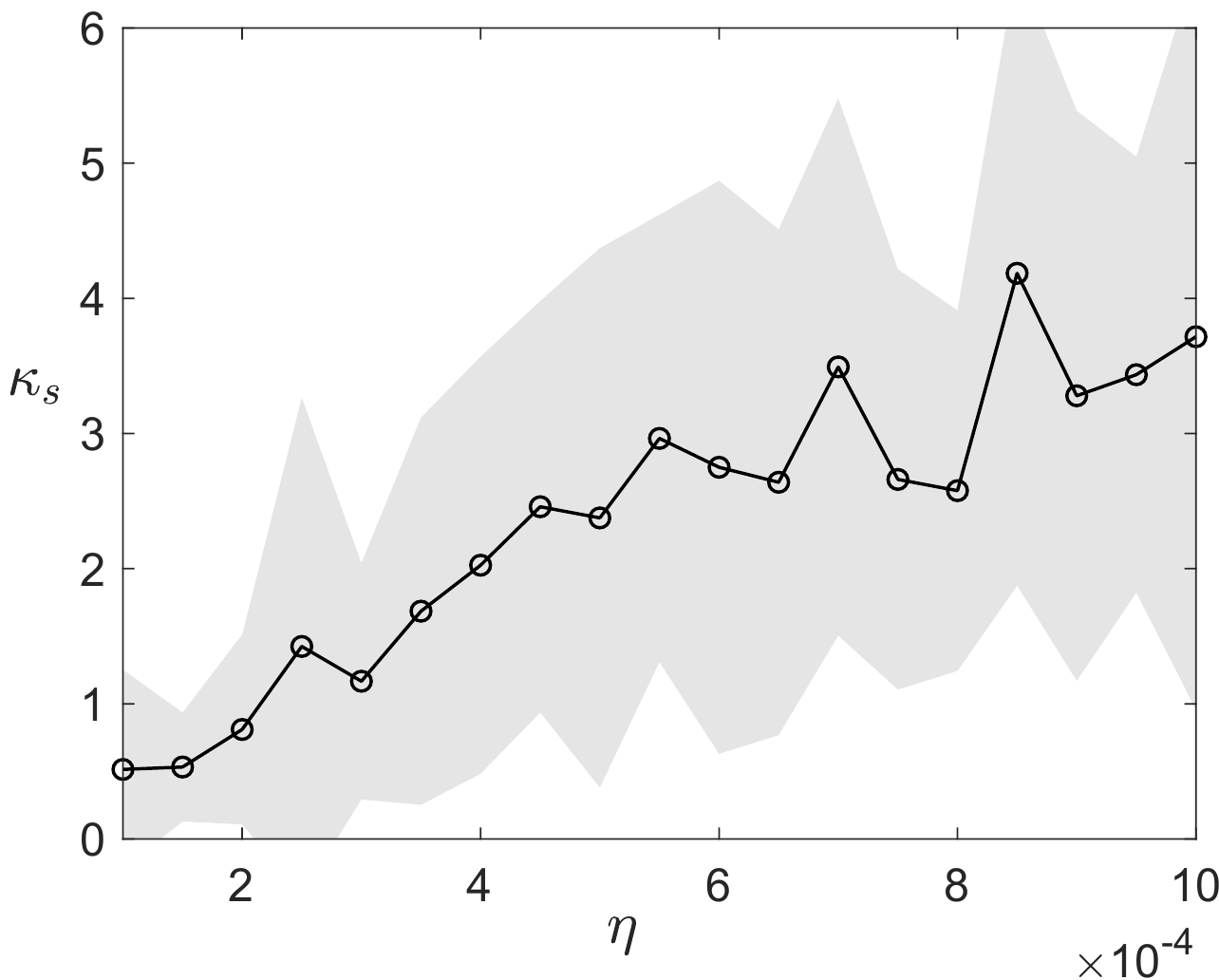}
\end{center}
\caption{Number of spurious connections $\kappa_s$ predicted by the LASSO in the presence of noise of intensity $\eta$, for a star network with $N=10$ and averaged over $20$ random initial conditions (shaded region represents the corresponding variance).} 
\label{fig:noise_case}
\end{figure}

 
  {We use as a measure of performance the number of recovered connections.  Suppose we have found equations of motion for the variables $\phi_i(t)$ in the form of a vector of coefficients ${w}^{(i)}$, $1\le i\le N$, where $w^{(i)}$ denotes the $i$th column of ${W}$. The norm of the function ${h}^{(i)}$, we can recover the strength of the coupling between node $i$ and the central node $1$ as}
\be 
\kappa_i=\| {h}^{(i)}\|_2=\sqrt{(c^{(i)}_{1i})^2+(d^{(i)}_{1i})^2}.
\ee 
 {We  define a quantity playing the role of effective total number of connections as }
\be 
\kappa=\frac{1}{\alpha} \sum_{i=2}^N \kappa_i ,
\ee
 {and the effective number of spurious connections (in general not an integer number),}
\be \kappa_s=\kappa-N.\ee
In Figure~\ref{fig:noise_case} we show how the performance of the method deteriorates as the amplitude of the noise increases, by plotting the effective number of spurious connections as a function of the noise intensity $\eta$, averaged over 50 random initial conditions (shaded region corresponds to standard deviation).

Equation~(\ref{phasenoise}) can be recast in the linear form  
\be \label{eq:linear_noisy_equation}
    v^{(i)} + \sqrt{2 \eta}~z^{(i)} = \Theta w^{(i)}, \quad i = 1, \dots, N,
\ee 
where $v^{(i)} \in \mathbb{R}^{n}$ corresponds to the Euler approximation of the time-derivative. Besides $z^{(i)} \in \mathbb{R}^{n}$ is a random variable whose each entry has the form of Equation (\ref{eq:Wiener_process}). 
Equation~(\ref{eq:linear_noisy_equation}) written in this form is similar to the noisy recovery case estimated by Cand\`es. Again, we take advantage of the relation between LASSO and the quadractically constrained basis pursuit problem.  

\begin{proposition}
Assume that $\delta_{2s} < \sqrt{2} - 1$. Given $\eta > 0$, $n > 0$ and $i = 1, \dots, N$, for each $\varepsilon > 0$ the following holds
\begin{align}
    \mathbb{P}\left(\|z^{(i)}\|_2 \leq \frac{\varepsilon}{\sqrt{2 \eta}}\right) \geq 1 - e^{-\frac{1}{4} \Big(\frac{\varepsilon}{\sqrt{2 \eta }}  - \sqrt{\frac{2}{\pi}} \sqrt{n} \Big)^2}
\end{align} 
So with high probability there exists the solution $w^{\star}$ of Equation (\ref{eq:quadractically_constrained_basis_pursuit}) satisfies
\[ 
\|w^* - w\|_2 \le C_0 s^{-1/2} \| w  - w_s\|_1 + C_1 \varepsilon
\] 
for constants $C_0$ and $C_1$

\end{proposition}

\begin{proof}
Fix a $\eta > 0$. We need to estimate how probable $\sqrt{2 \eta}z^{(i)}$ is $L_2$ bounded by a constant $\varepsilon$.  We drop the dependence of $(i)$. Note that 
\begin{linenomath*} \be
    \sqrt{2 \eta} \|z\|_2 \leq \sqrt{2 \eta} (\|\xi\|_2 + \|\zeta\|_2)
 \ee \end{linenomath*}
where each $\xi$ and $\zeta$ are Gaussian random vectors. So, we can estimate the expected value of this $L_2$ norm: $\mathbb{E}(\|\xi\|_2) \geq  \sqrt{\frac{2}{\pi} n}$~\cite[Proposition 8.1]{Foucart_mathematical_compressive_sensing} and the same for $\zeta$. 
We use the concentration of measure for Gaussian random vector~\cite[Theorem 8.34]{Foucart_mathematical_compressive_sensing}. Since the norm $L_2$ is a Lipschitz function with constant $1$, the estimate follows. Then $\|\xi\|_2 \leq \frac{\varepsilon}{2\sqrt{2 \eta}}$ holds with probability
\begin{linenomath*} \be
    1 - \mbox{exp} \left(- \frac{1}{4}\left(\frac{\varepsilon}{\sqrt{2 \eta }}  - \sqrt{\frac{2}{\pi}} \sqrt{n}\right)^2 \right).
 \ee \end{linenomath*}
Consequently, we can apply Cand\`es' estimate and the statement is proved. 
\end{proof}

\appendix

{  

\section{Proof of Proposition \ref{prop:L_2_instability}}\label{Insta}

We first need \incorrect{two} preliminary results

\begin{OBS}\label{lemma:svd_not_kernel}
Given $S \in \mathbb{R}^{s_1 \times s_2}$ with $1 \leq rank(S) \leq \min\{s_1, s_2\}$, then 
\begin{linenomath*} \be
    \|S x\|_2 \geq \sigma_{min}(S) \|x\|_2, \forall~x \in \mathbb{R}^{s_2}\backslash \{\ker S\}, \nonumber
    \ee
\end{linenomath*} 
where $\sigma_{\min}(S)$ is the minimum singular value of $S$~\cite[Fact 9.13.1]{bernstein2009matrix}.
\end{OBS}

\begin{OBS}\label{lemma:generic_z}
 Let $A \in \mathbb{R}^{n \times p}$ and $B \in \mathbb{R}^{n \times q}$ be column full rank matrices with $n > p + q$, and $\varepsilon \in (0, 1)$. Let $r = \min\{p, q\}$. We have the following:
 \begin{enumerate}
     \item[i)] For generic $z \in \mathbb{R}^{n}$: $ z \not \in (\im~B)^{\perp}$.  
 The map $H(z) = B^{\dagger} z$ is not constant thus Leb$(H^{-1}(0)) = 0$ and $\mathbb{R}^{n}\backslash H^{-1}(0)$ is a generic set. \par 
     \item[ii)] For generic $z \in (\im~A)^{\perp}$ with $\|z\|_2 = \varepsilon$, there exists $K(B, \varepsilon) > 0$ such that  
\[
\|B^{\dagger} z\|_2 \geq K \cos(\beta_r), 
\]
where $\beta_r$ is the largest principle angle between $(\im~A)^{\perp}$ and $\im~B$. Indeed, let $Q_B R_B$ be the $QR$ decomposition of the matrix $B$ and $Q_{A^{\perp}}$ an orthonormal matrix whose columns form an orthonormal basis to $(\im~A)^{\perp}$. Hence, there exists a unique $v \not = 0$ such that $z = Q_{A^{\perp}} v$ and $\|v\|_2 = \varepsilon$. Applying Remark~\ref{lemma:svd_not_kernel} and previous item (i), generically, we have 
 \begin{linenomath*} \be
     \|B^{\dagger} z\|_2 \geq \sigma_{\min}(R_B) \sigma_{\min} (Q_B^{\dagger}Q_{A^{\perp}}) \|v\|_2.
  \ee \end{linenomath*}
 By \cite[Theorem 2.1]{PABS} the principal angles $\beta$'s 
 between subspaces $(\im~A)^{\perp}$ and $\im~B$ are 
 \begin{linenomath*} 
 \begin{equation}
(\cos(\beta_1),  \dots  ,\cos(\beta_r)) = 
(     \sigma_{\max}(Q_B^{\dagger} Q_{A^{\perp}}),  \dots , \sigma_{\min}(Q_B^{\dagger} Q_{A^{\perp}})),
 \end{equation}
 \end{linenomath*}
 where $\beta_k \in [0, \frac{\pi}{2}]$, $k = 1, \dots, r$ and $\beta_{k} < \beta_{k + 1}$, $k = 1, \dots, r - 1$. In particular, the cosine of the largest angle between $(\im~A)^{\perp}$ and $\im~B$ is given as follows
\begin{linenomath*} \be
\cos \beta_r = \sigma_{\min} (Q_{B}^{\dagger} Q_{A^{\perp}}).
 \ee \end{linenomath*}
So, there exists $K(B, \varepsilon) > 0$ such that
\begin{linenomath*} \be
     \|B^{\dagger} z\|_2 \geq K \cos(\beta_r).
  \ee \end{linenomath*}
 \end{enumerate}
\end{OBS}
\begin{lemma}\label{lemma:properties_of_M}
Let $A \in \mathbb{R}^{n \times p}$ and $B \in \mathbb{R}^{n \times q}$ be column full rank matrices with $n > p + q$. Let $r = \min\{p, q\}$ then $\beta_r \not = \frac{\pi}{2}$ be the largest principle angle between the subspaces $(\im~A)^{\perp}$ and $\im~B$. Consider
\begin{linenomath*} \be
M = B^{\dagger}B - B^{\dagger} A (A^{\dagger} A)^{-1} A^{\dagger} B
\ee
\end{linenomath*} 
then there exists a constant $K > 0$ such that 
\begin{linenomath*} \be
    \sigma_{\min}\big(M^{-1}\big) \geq \frac{1}{K \cos^2(\beta_r)}. 
\ee
\end{linenomath*} 
\end{lemma}

\begin{proof}
Let $E$ and $F$ be the orthogonal projection onto the $\im~A$ and $(\im~A)^{\perp}$, respectively. So, using that $A$ is column full rank we can write $M$ as follows:
\begin{linenomath*} 
\begin{eqnarray*}
    M &= B^{\dagger}B - B^{\dagger} A (A^{\dagger} A)^{-1} A^{\dagger} B \\
    &= B^{\dagger} (\mathbf{I} - A A^{+})B\\
    &= B^{\dagger} FB \in \mathbb{R}^{q \times q}.
 \end{eqnarray*} 
 \end{linenomath*}
The orthogonal projections have the following formulas: $E = Q_{A}Q_A^{\dagger}$ and $F = Q_{A^{\perp}}Q_{A^{\perp}}^{\dagger}$, where $Q_{A}$ and $Q_{A^{\perp}}$ are orthonormal matrices whose columns form an orthonormal basis to $\im~A$ and $(\im~A)^{\perp}$, respectively. Besides, let us denote $B = Q_{B} R_B$ the QR decomposition of the matrix $B$. Using this notation and inequality of the singular value~\cite[Theorem 3.3.14]{horn_2011} we can split up the maximum singular value as follows:
\begin{linenomath*} 
\be
\label{eq:inequality_M_}
     \sigma_{\max}(B^{\dagger} FB) \leq \sigma_{\max}(R_B^{\dagger}) \sigma_{\max}(Q_{B}^{\dagger}Q_{A^{\perp}})\sigma_{\max} (Q_{A^{\perp}}^{\dagger} Q_B) \sigma_{\max}(R_B).
\ee
\end{linenomath*} 
Again by \cite[Theorem 2.1]{PABS}  the least angle between $(\im~A)^{\perp}$ and $\im~B$ is given as follows
\begin{linenomath*} 
\be
\label{eq:PABS}
\cos \beta_1 = \sigma_{\max} (Q_{B}^{\dagger} Q_{A^{\perp}}) = \sigma_{\max} (Q_{A^{\perp}}^{\dagger} Q_B),
\ee
\end{linenomath*} 
 Also $\beta_1 < \beta_r $ so $\beta_1 \not = \frac{\pi}{2}$ and Equation~\eqref{eq:PABS} is not zero. If we use that $\sigma_{\max}(B^{\dagger} FB) = \frac{1}{\sigma_{\min}((B^{\dagger} FB)^{-1})}$ and Equation~\eqref{eq:PABS} into Equation~\eqref{eq:inequality_M_} we obtain
\begin{linenomath*} \be
    \sigma_{\min}((B^{\dagger} FB)^{-1}) \geq \frac{1}{K\cos^{2}(\beta_1)},
 \ee \end{linenomath*}
where we use that $\sigma_{\max}(R_B)$ can be bounded by a constant $K > 0$. Since $\cos$ is decreasing in the interval $[0, \frac{\pi}{2}]$, we can replace $\beta_1$ by the largest principle angle $\beta_r$, and the claim follows.
\end{proof}

\begin{proof}[Proof of Proposition~\ref{prop:L_2_instability}] 
Note that since $C$ is column full rank, we can write the solution of Equation~\eqref{eq:partitioned_l2_minization} as
\begin{linenomath*} 
\be
\begin{pmatrix}\hat w_1 \\ \hat w_2\end{pmatrix} = C^{+} (b + z) 
= \begin{pmatrix} A^{\dagger} A & A^{\dagger} B \\
                    B^{\dagger} A & B^{\dagger}B \end{pmatrix}^{-1} \begin{pmatrix} A^{\dagger} b \\ 
                    B^{\dagger} b + B^{\dagger} z\end{pmatrix},
 \ee \end{linenomath*}
where we use that $z \in (\im~ A)^{\perp}$ implies that $A^{\dagger} z = 0$.  
Using the analytic inversion formula~\cite{bernstein2009matrix}, we obtain
\begin{linenomath*} 
\be
(C^{\dagger} C)^{-1} =  
 \begin{pmatrix} (A^{\dagger} A)^{-1} + (A^{\dagger} A)^{-1} (A^{\dagger}B) M^{-1} (B^{\dagger} A) (A^{\dagger} A)^{-1} & - (A^{\dagger} A)^{-1} A^{\dagger} B M^{-1} \\
                    -M^{-1} (B^{\dagger}A)(A^{\dagger} A)^{-1} & M^{-1} \end{pmatrix}
 \ee \end{linenomath*}
where $M^{-1} = (B^{\dagger}B - B^{\dagger} A (A^{\dagger} A)^{-1} A^{\dagger} B)^{-1}$. 
Since $A$ and $B$ are column full rank, we can use the formula of $A^{+}$ and $A A^{+}$ is a projector onto the $\im~A$. So, we obtain
\begin{linenomath*} 
\be
\begin{pmatrix}\hat w_1 \\ \hat w_2\end{pmatrix}  = C^{+} (b + z) \\
= \begin{pmatrix} A^{+} b - A^{+} B M^{-1} B^{\dagger}z \\
                 M^{-1} B^{\dagger}z\end{pmatrix}
 \ee \end{linenomath*}
where we used that $b \in \im~A$. \par 
We aim at calculating how much the solution $x^{*}$ is perturbed, so 
\begin{linenomath*} \be
    \| x^* - \hat w_1 \|_2 = \|  A^{+} B M^{-1} B^{\dagger} z \|_2.
 \ee \end{linenomath*}
Since $\beta_1 > 0$ we have $B M^{-1} B^{\dagger} z \not \in (\im~A)^{\perp}$, so using Remark~\ref{lemma:svd_not_kernel} for $A^{+}$ we obtain
\begin{linenomath*} \be
    \| x^* - \hat w_1 \|_2 \geq  \sigma_{\min}(A^{+}) \|B M^{-1} B^{\dagger} z \|_2.
 \ee \end{linenomath*}
By item (i) of Remark~\ref{lemma:generic_z} for generic $z \not = 0$ we have $B^{\dagger} z \not = 0$.
Recall that $M^{-1} B^{\dagger} z \in (\ker B)^{\perp}$. Thus, Remark~\ref{lemma:svd_not_kernel} is valid for $B$ and $M^{-1}$ and we obtain 
\begin{linenomath*} 
\be
    \| x^* - \hat w_1 \|_2 \geq \sigma_{\min}(A^{+}) \sigma_{\min}(B)  \sigma_{\min}(M^{-1}) \| B^{\dagger} z\|_2. \nonumber
\ee
\end{linenomath*} 
Observe that $\sigma_{\min}(A^{+}), \sigma_{\min}(B)  > 0$ from the full rank condition on the matrices $A$ and $B$, so there exists $K_1 > 0$ given by 
\begin{linenomath*} 
\be
K_1 = \min\{\sigma_{\min}(A^{+}), \sigma_{\min}(B)\}.
\ee
\end{linenomath*} 
Moreover, by item (ii) of Lemma~\ref{lemma:generic_z} we have: $\|B^{\dagger}z\|_2 \geq K_2 \cos(\beta_r)$. Hence, by \cite[Property 2.1]{PABS} if $|\beta_r - \pi/2| < \varepsilon$ for sufficiently small $\varepsilon$, using Lemma~\ref{lemma:properties_of_M} there exists $K> 0$ we obtain
\begin{linenomath*} 
\be
    \| x^* - \hat w_1 \|_2 \geq \frac{K_1 K_2 }{K \cos \beta_r} > 0
\ee
\end{linenomath*} 
and the statement holds.
\end{proof}
}

\bibliographystyle{unsrt}

\end{document}